\documentclass[final,reqno,oneeqnum,onefignum,onetabnum]{siamart190516}

\usepackage{amssymb}
\usepackage{mathtools}
\usepackage{booktabs}
\usepackage{multirow}
\usepackage{microtype}
\usepackage{comment}
\usepackage{centernot}
\usepackage{cancel}

\usepackage{tikz}
\usetikzlibrary{decorations.markings}

 \usepackage{pgfplots}
 \usepgfplotslibrary{polar}
 \pgfplotsset{compat=1.14}

\usepackage[T1]{fontenc}
\usepackage{lmodern}

\usepackage[round,sort]{natbib}
\usepackage{hyperref}

\title
{Fast high-dimensional integration \\using tensor networks}

\headers
{Fast high-dimensional integration using tensor networks} 
{S.~Cassel}

\author{Sebastian~Cassel\thanks{BNP Paribas, 10 Harewood Avenue, London, NW1 6AA, UK 
		(\email{sebastian.cassel@protonmail.com}).
}}

\begin{document}
\maketitle
\thispagestyle{empty} 

\begin{abstract}
The design and application of regression-free tensor network representations for integration is presented. Tensor network methods are demonstrated to outperform Monte Carlo for test problems, and exponential convergence is shown to be achievable for a non-analytic integrand.
\end{abstract}






\noindent 
A common task in scientific disciplines is to calculate high-dimensional integrals, for~example when solving integral/differential equations or evaluating expectations over probability distributions. Such problems are extensive in finance for derivative pricing and risk analysis, where dependencies on large numbers of state variables and time horizons can exist. Furthermore, these problems routinely demand significant computational resources and so techniques for improving efficiency generally offer notable advantages. 

Historically, Monte Carlo has been the leading method for high-dimensional numerical integration. Although Monte Carlo errors decay as $\mathcal{O}[n^{-\frac12}]$ irrespective of the dimension (given $n$ samples), the convergence rate is very slow: to suppress the error by~$10^{m}$, the sample count and thus runtime must multiply by~$10^{2m}$. A future opportunity is to use quantum amplitude estimation \citep{montanaro2015quantum} on quantum computers to achieve linear convergence $\mathcal{O}[n^{-1}]$ for general integrals. However, such convergence rates are not necessarily optimal for a given problem.

This article emphasises the usefulness of tensor networks for high-dimensional integration on classical computers, demonstrating convergence rates that outperform Monte Carlo and that exponential convergence $\mathcal{O}[e^{-\alpha n}]$ is even achievable for certain integrands.


\vspace{4mm}\noindent
\textbf{Tensor network introduction}

\begin{figure}[!htbp]
	\centering
	\begin{tikzpicture}
		\def \n {6}
		\def \radius {0.6cm}
		\def \rk {1.7cm}
		\def \margin {16} 
		\node[draw, fill=white, scale=0.5, fill=gray!20, circle] (node1) at ({90 - 0*60}:\radius) {};
		\node[draw, fill=white, scale=0.5, fill=gray!20, circle] (node2) at ({90 - 1*60}:\radius) {};
		\node[draw, fill=white, scale=0.5, fill=gray!20, circle] (node3) at ({90 - 2*60}:\radius) {};
		\node[draw, fill=white, scale=0.5, fill=gray!20, circle] (node4) at ({90 - 3*60}:\radius) {};
		\node[draw, fill=white, scale=0.5, fill=gray!20, circle] (node5) at ({90 - 4*60}:\radius) {};
		\node[draw, fill=white, scale=0.5, fill=gray!20, circle] (node6) at ({90 - 5*60}:\radius) {};
		\node[fill=white, inner sep=0.5mm, circle] (k1) at ({90 - 0*60}:\rk) {$k_1^{}$};
		\node[fill=white, inner sep=0.5mm, circle] (k2) at ({90 - 1*60}:\rk) {$k_2^{}$};
		\node[fill=white, inner sep=0.5mm, circle] (k3) at ({90 - 2*60}:\rk) {$k_3^{}$};
		\node[fill=white, inner sep=0.5mm, circle] (k4) at ({90 - 3*60}:\rk) {$k_4^{}$};
		\node[fill=white, inner sep=0.5mm, circle] (k5) at ({90 - 4*60}:\rk) {$k_5^{}$};
		\node[fill=white, inner sep=0.5mm, circle] (k6) at ({90 - 5*60}:\rk) {$k_6^{}$};
		\draw[-] (node1) -- (k1);
		\draw[-] (node2) -- (k2);
		\draw[-] (node3) -- (k3);
		\draw[-] (node4) -- (k4);
		\draw[-] (node5) -- (k5);
		\draw[-] (node6) -- (k6);
		\draw[fill=gray!5] (0,0) circle (0.9cm);
		\node[circle, scale=1.0] (core) at (0:0) {$T_{ k_1^{} \cdots k_6^{} }^{}$};
		\node[circle, scale=1.5] at (0:3.0cm) {$\Longrightarrow$};
	\end{tikzpicture}
	\hspace{5mm}
	\begin{tikzpicture}
		\def \n {6}
		\def \radius {0.6cm}
		\def \rk {1.7cm}
		\def \margin {16} 
		\draw[fill=gray!5, dashed] (0,0) circle (0.9cm);
		%
		\node[draw, fill=white, scale=0.5, fill=gray!20, circle] (node1) at ({90 - 0*60}:\radius) {};
		\node[draw, fill=white, scale=0.5, fill=gray!20, circle] (node2) at ({90 - 1*60}:\radius) {};
		\node[draw, fill=white, scale=0.5, fill=gray!20, circle] (node3) at ({90 - 2*60}:\radius) {};
		\node[draw, fill=white, scale=0.5, fill=gray!20, circle] (node4) at ({90 - 3*60}:\radius) {};
		\node[draw, fill=white, scale=0.5, fill=gray!20, circle] (node5) at ({90 - 4*60}:\radius) {};
		\node[draw, fill=white, scale=0.5, fill=gray!20, circle] (node6) at ({90 - 5*60}:\radius) {};
		\node[fill=white, inner sep=0.5mm, circle] (k1) at ({90 - 0*60}:\rk) {$k_1^{}$};
		\node[fill=white, inner sep=0.5mm, circle] (k2) at ({90 - 1*60}:\rk) {$k_2^{}$};
		\node[fill=white, inner sep=0.5mm, circle] (k3) at ({90 - 2*60}:\rk) {$k_3^{}$};
		\node[fill=white, inner sep=0.5mm, circle] (k4) at ({90 - 3*60}:\rk) {$k_4^{}$};
		\node[fill=white, inner sep=0.5mm, circle] (k5) at ({90 - 4*60}:\rk) {$k_5^{}$};
		\node[fill=white, inner sep=0.5mm, circle] (k6) at ({90 - 5*60}:\rk) {$k_6^{}$};
		\draw[-] (node1) -- (node2) node[draw, fill=white, scale=0.5, circle, pos=0.5] {};
		\draw[-] (node2) -- (node3) node[draw, fill=white, scale=0.5, circle, pos=0.5] {};
		\draw[-] (node3) -- (node4) node[draw, fill=white, scale=0.5, circle, pos=0.5] {};
		\draw[-] (node4) -- (node5) node[draw, fill=white, scale=0.5, circle, pos=0.5] {};
		\draw[-] (node5) -- (node6) node[draw, fill=white, scale=0.5, circle, pos=0.5] {};
		\draw[-] (node1) -- (k1);
		\draw[-] (node2) -- (k2);
		\draw[-] (node3) -- (k3);
		\draw[-] (node4) -- (k4);
		\draw[-] (node5) -- (k5);
		\draw[-] (node6) -- (k6);
	\end{tikzpicture}
	\\[3mm]
	{ \raggedright 
		\hspace{65mm}
		$\displaystyle \sum_{\{ \ell_i^{}, m_i^{} \}}^{} 
		\tilde{T}_{k_1^{} \ell_1^{}}^{(a)} \, \tilde{T}_{\ell_1^{} m_1^{}}^{(b)} \, \tilde{T}_{m_1^{} k_2^{} \ell_2^{}}^{(c)} \, \cdots$ }
	\caption{Graphical representation of an example tensor network}
	\vspace{-2mm}
	\label{fig:network}
\end{figure}
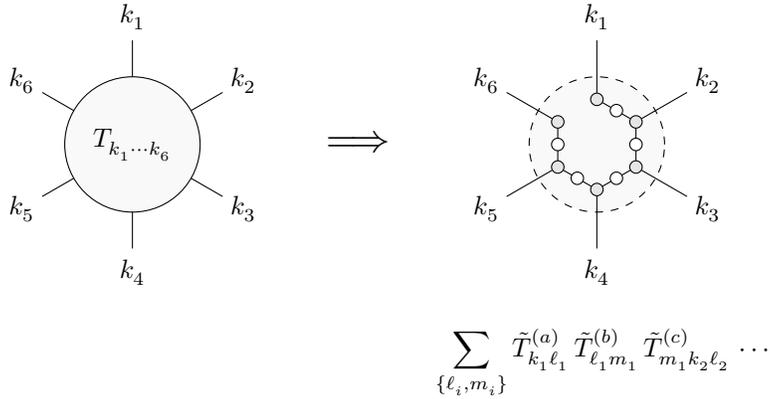

\noindent
The structure of a tensor network representation is related to an (arbitrary) integral representation of a given multivariate function:
\begin{eqnarray} 
	f [\mathbf{x}] &~=~& \int T [\mathbf{k}] \, \Bigg( \prod_{i}^{d} \Gamma^{(i)} [k_i^{}, x_i^{}]  \Bigg) \, d\mathbf{k}
\end{eqnarray}
On discretising the integral, a core tensor $T$ is formed which generally has an exponentially large number of components (since if each $k_i^{}$ index assumes $n$ values there are $n^d$ components).
\begin{equation}
	f [\mathbf{x}]  ~\approx \sum_{ \text{T-indices} }^{} T_{k_1^{} \cdots k_d^{}}^{} \, \Bigg( \prod_{i}^{d} \Gamma_{k_i^{}}^{(i)} [x_i^{}] \Bigg)
	\label{eq:core} 
\end{equation} 
In order to then control the complexity, a network of low rank tensors can be chosen to replace the core tensor. For example, \cref{fig:network} illustrates a particular network, where each node corresponds to a tensor and the number of respective tensor indices corresponds to the number of connecting lines. 
Although an infinite sum over network indices may be needed to fully reproduce the continuum limit of the core tensor, representations with rapidly decaying approximation errors are sought for practical application.

The separable form of \cref{eq:core} is strikingly convenient for integration over $\{x_i^{}\}$, as one-dimension integration techniques can be applied. However, it is essential that low-complexity tensor network representations can be efficiently formed for the approach to be useful. In the fields of quantum physics and chemistry, tensor network regressions have been widely used over past decades. In finance, tensor network regressions have recently been applied by \citet{glau2020low} and \citet{antonov2021alternatives}. A limitation for regression though is that its calibration tends to be computationally intensive, restricting the network size. This article highlights regression-free methods for forming tensor network representations, for which larger networks can be supported. 

\vspace{4mm}\noindent
\textbf{Tensor Train Cross (TT-X) network}

\noindent
Constructive methods for forming tensor networks can be based on interpolation, noting that interpolation constraints for multivariate functions may be lines (or even hypersurfaces) instead of simply points. For example, in two dimensions, a function that interpolates with respect to co-ordinate lines intersecting at node $\mathbf{s}_\ast^{}$ is given by:
\begin{eqnarray}
	f [x_1^{} , x_2^{}] &~\approx~& \frac{ f[x_1^{}, s_{\ast 2}^{}] ~f[s_{\ast 1}^{}, x_2^{}] }{f[s_{\ast 1}^{}, s_{\ast 2}^{}]} 
\end{eqnarray}
More generally, a two-dimensional function with co-ordinate line constraints intersecting at each node $\{ \mathbf{s}_i^{} \}$ is given by the following expression:
\begin{align}
	f [x_1^{} , x_2^{}] ~\approx~ \sum_{i,j}^{} f [x_1^{}, s_{i 2}^{}] ~ ( \mathbf{Q}_{1,2}^{-1} )_{ij}^{} ~  f [s_{j 1}^{}, x_2^{}] 
	\label{eq:ttx2}
	\\[1mm]
	\text{where} \hspace{2mm}
	( \mathbf{Q}_{1,2}^{} )_{k\ell}^{} \,=\,  f [s_{k 1}^{}, s_{\ell 2}^{}]  \hspace{10mm} \phantom{,} \nonumber
\end{align}
Although interpolation properties may break if nodes are selected such that~$\mathbf{Q}$ is singular, an approximation can still be formed by applying a pseudo-inverse of~$\mathbf{Q}$ in~\cref{eq:ttx2}.

On moving to three dimensions, an interpolating approximation takes the form:
\begin{eqnarray} 
	f [x_1^{} , x_2^{}, x_3^{}] &~\approx& \sum_{\text{T-indices}}^{} T_{ijk\ell mn}^{} ~ f [x_1^{}, s_{i 2}^{}, s_{j 3}^{}] ~  f [s_{k 1}^{}, x_2^{}, s_{\ell 3}^{}] ~ f [s_{m 1}^{}, s_{n 2}^{}, x_3^{}] 
\end{eqnarray}
The respective core tensor $T$ attains exponential complexity in higher dimensions, and so tensor networks then become useful to control the complexity. Such an approach is equivalent to controlling the number of interpolation constraints applied. 

\begin{figure}[!htbp]
	\centering
	\begin{tikzpicture}
		\pgfplotsset{axis line style={gray!70} }
		\begin{axis}[xmin=0, xmax=1, ymin=0, ymax=1, ticks=none, width=5.5cm, height=5.5cm,
			axis background/.style={fill=gray!10}, name = plot2d] 
			\addplot[gray] coordinates { (0.2, 0) (0.2, 1) };
			\addplot[gray] coordinates { (0.33, 0) (0.33, 1) };
			\addplot[gray] coordinates { (0.75, 0) (0.75, 1) };
			\addplot[gray] coordinates { (0, 0.40) (1, 0.40) };
			\addplot[gray] coordinates { (0, 0.25) (1, 0.25) };
			\addplot[gray] coordinates { (0, 0.70) (1, 0.70) };
			\addplot [only marks,mark=*, mark options={scale=1.2}] coordinates { (0.2,0.40) (0.33,0.25) (0.75,0.70) };
			\addplot [only marks,mark=*, mark options={scale=1.2, fill=white}] coordinates 
			{ (0.2,0.25) (0.2,0.70) (0.33,0.40) (0.33,0.70) (0.75,0.25) (0.75,0.40) };
			\node[label={$\mathbf{s}_1^{}$}] at (0.12,0.36){};
			\node[label={$\mathbf{s}_2^{}$}] at (0.41,0.07){};
			\node[label={$\mathbf{s}_3^{}$}] at (0.83,0.66){};
		\end{axis}
		\begin{axis}[clip=false, ticks=none, width=6cm, height=6cm,
			axis background/.style={fill=gray!10}, 
			name=plot3d, at=(plot2d.right of east), anchor=left of west, zlabel={\phantom{,}}, zlabel shift=15mm]
			\addplot3[mark=none, black] coordinates {(0.1,0,0.3) (0.1,1,0.3)};
			\addplot3[mark=none, black] coordinates {(0,0.2,0.3) (1,0.2,0.3)};
			\addplot3[mark=none, black] coordinates {(0.1,0.2,0) (0.1,0.2,1)};
			\addplot3[mark=none, black] coordinates {(0.8,0,0.85) (0.8,1,0.85)};
			\addplot3[mark=none, black] coordinates {(0,0.75,0.85) (1,0.75,0.85)};
			\addplot3[mark=none, black] coordinates {(0.8,0.75,0) (0.8,0.75,1)};
			\addplot3[only marks,mark=*, black, mark options={scale=1.0}] coordinates { (0.1,0.2,0.3) (0.8,0.75,0.85) };
			\addplot3[only marks,mark=*, black, mark options={scale=1.0, fill=white}] coordinates 
			{ (0.1,0.75,0.85) (0.8,0.2,0.3) (0.1,0.2,0.85) (0.8,0.75,0.3)   };
		\end{axis}
	\end{tikzpicture}
	\caption{Line interpolation constraints within TT-X representations in two and three dimensions. The dark points correspond to nodes $\{ \mathbf{s}_i^{} \}$, and light points to off-diagonal elements in the $\mathbf{Q}$ matrices.}
	\label{fig:lines}
\end{figure}
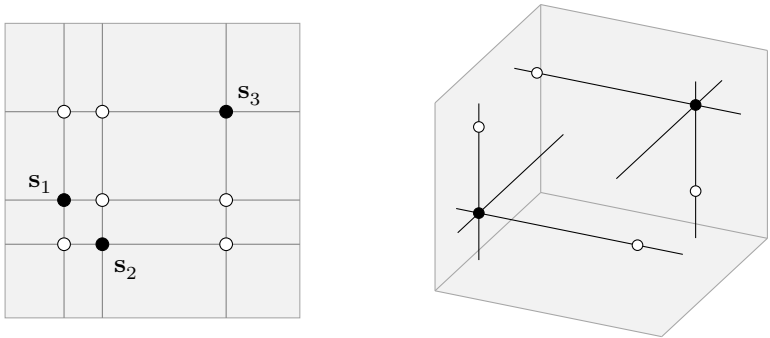

For an arbitrary number of dimensions, the minimal-complexity form supporting interpolation with respect to lines is given by the tensor train cross representation \citep{oseledets2010tt}:
\begin{equation}
	f [x_1^{}, \,\cdots, x_d^{}] ~\approx~ 
	\mathbf{F}_1^{} [x_1^{}] ~ \mathbf{Q}_{1,2}^{-1} ~ \mathbf{F}_2^{} [x_2^{}]  
	~\cdots~ \mathbf{Q}_{d-1,d}^{-1}~
	\mathbf{F}_d^{} [x_d^{}] 
	\label{eq:ttx}
\end{equation}
\vspace{-6mm}
\begin{eqnarray}
	( \mathbf{F}_a^{} [x_a^{}] )_{k\ell}^{} &~=~&
	f[ s_{k \lhd}^{} , x_a^{} , s_{\ell \, \rhd}^{} ]
	\\[2mm]
	( \mathbf{Q}_{a,b}^{} )_{k\ell}^{} &~=~&
	f[ s_{k \lhd}^{} , s_{ka}^{} , s_{\ell \, \rhd}^{} ]
\end{eqnarray}
where the notation introduced is defined as follows (for free variable $z_a^{}$):
\begin{eqnarray}
	f[ s_{k \lhd}^{} , z_a^{} , s_{\ell \, \rhd}^{} ]
	&~=~&
	f[s_{k1}^{}, \ldots, s_{k (a-1)}^{}, \, z_a^{}, \, s_{\ell (a+1)}^{}, \ldots, s_{\ell d}^{}]  
\end{eqnarray}
The tensor train cross representation in \cref{eq:ttx} is composed of matrix functions and matrix connections associated with a node set $\{ \mathbf{s}_i^{} \}$. Also to note, the edge matrix functions $\mathbf{F}_1^{}$ and $\mathbf{F}_d^{}$ necessarily only have one row and one column respectively so that a scalar function is formed by \cref{eq:ttx}:
\begin{eqnarray}
	f[ s_{k \lhd}^{} , z_1^{} , s_{\ell \, \rhd}^{} ]
	&~=~&
	f[  z_1^{}, \, s_{\ell 2}^{}, \ldots, s_{\ell d}^{}]  \\[2mm]
	f[ s_{k \lhd}^{} , z_d^{} , s_{\ell \, \rhd}^{} ]
	&~=~& f[s_{k1}^{}, \ldots, s_{k (d-1)}^{}, z_d^{} ]  
\end{eqnarray}
As illustrated by~\cref{fig:lines}, when moving from two to three (higher) dimensions the constraint intersections (`crossings') at nodes are maintained, but off-node intersections are generally lost.

In forming the tensor train cross representation, the following choices are generally available to influence the quality of the approximation:
\vspace{1mm}
\begin{itemize}
	\item choice of co-ordinate system
	\item choice of dimension ordering
	\item choice of node count and placement
\end{itemize}
\vspace{1mm}
Focusing on node placement, random selection is generally disfavoured as approximation errors are then uncontrolled. In the paper by \citet{oseledets2010tt}, it is proposed to select nodes that maximise a measure related to determinants of QR matrix decompositions of the representation terms. Such an approach is motivated to find dominant modes of the representation, but it attracts a significant computational burden. For the results in this article, successive nodes are simply chosen to be in the neighbourhood of large approximation errors. A local search can be performed to improve the error suppression, but it is not necessary to do so.

For practical use, the low complexity of a tensor train cross representation is a key advantage. Given $n$ nodes, there are $\mathcal{O} [d n^2 ]$ components that need to be calculated: the scaling is linear (not exponential) in dimension, and quadratic in node count. Importantly, it is also a constructive method so no calibration is needed. However, some care is still necessary for the node placement.

\vspace{4mm}\noindent
\textbf{Series acceleration: Aitken extrapolation}

\noindent
As tensor networks may be constructed deterministically, it can be possible to apply series acceleration techniques to the sequence of results in order to deduce asymptotic limits faster. In the following section, Aitken extrapolation will be applied to the TT-X integration results, and so this technique is now briefly described.

For a sequence of $\{ \psi_i^{} \}$ values, a telescoping relationship can be formed:
\begin{eqnarray}
	\psi_{n} &~=~&
	\psi_0^{} + \sum_{i=0}^{n-1} \, g_i^{} \hspace{5mm} \text{where} \hspace{2mm} g_i^{} \,=\,
	\left( \psi_{i+1}^{} - \psi_{i}^{} \right)
\end{eqnarray}
If $g_i^{}$ is interpreted as a function value $g [ \psi_i^{} ]$, such a sequence is equivalent to fixed-point iterations. The asymptotic value $\psi_{\infty}^{}$ then corresponds to where $g [ \psi_\infty^{} ]$ is zero, and so the secant method can be applied:
\begin{eqnarray}
	\psi_\infty^{} &~\approx~&
	\psi_i^{} - \frac{ g_i^{} }{ g_i^{\prime}  }
	\hspace{6mm} \text{where} \hspace{2mm} g_i^\prime \,\approx~ \frac{ g_{i}^{} - g_{i-1}^{} }{ \psi_{i}^{} - \psi_{i-1}^{} }
\end{eqnarray}
The above formula defines Aitken extrapolation. The suitability of a given series acceleration technique will depend on the integration problem and TT-X node selection procedure, but such methods can not be applied to stochastic estimates as produced by Monte Carlo (or quantum amplitude estimation).

\vspace{4mm}\noindent
\textbf{Basket option valuation}

\noindent
The TT-X representation can be adopted for any function, and its effectiveness is now demonstrated for European basket option valuation. Assuming Black-Scholes dynamics, the option value $\psi$ is given by the solution of a linear differential equation. Such solutions can be expressed as an integral of a Green's function $G$ and relevant source function $\phi$ as follows:
\begin{align}
	\psi [\mathbf{x}, t] &~=~ \int G [\mathbf{x}, t ; \mathbf{x}^{\prime}, t^\prime] \, \phi [\mathbf{x}^{\prime},t^\prime] \, d\mathbf{x}^\prime \, dt^\prime
	\label{eq:psi}
	\\[4mm]
	\phi [\mathbf{x}^\prime, t^\prime] &~=~ \, \delta (t^\prime - t_{\ast}^{} ) \,\,
	\max \Bigg[ \, 0, \Bigg( \sum_{i}^{d} \omega_i^{} \, e^{x_i^{\prime}} \Bigg) - K \, \Bigg]
	\label{eq:phi}
	\\[5mm]
	G [\mathbf{x}, t ; \mathbf{x}^{\prime}, t^\prime] &~=~
	\frac{
	\exp \left[ - \frac12\, (\mathbf{x}^\prime - \boldsymbol{\mu})^T\, \boldsymbol{\Sigma}^{-1} (\mathbf{x}^\prime - \boldsymbol{\mu}) 
	- r\, (t^\prime - t) 
	\right]
	}{
		\sqrt{ \det 2\pi \boldsymbol{\Sigma} }}
	\label{eq:G}
\end{align}
where $\boldsymbol{\mu}$ and $\boldsymbol{\Sigma}$ generally depend on $\{ \mathbf{x}, t , t^\prime \}$ and model parameters. In the following tests, the default parameter settings listed in \cref{tbl:parameters} are applied 
unless $\rho$ or $d$ is specifically referenced.

\begin{table}[!htbp]
	\centering
	\begin{tabular}{llcc}
		\toprule
		& Parameter & Default & Alternative \\ \midrule
		$d$ & Basket dimension & 10 & 100 \\
		$\omega_i^{}$ & Basket weight & $1/d$ \\
		$\mu_i^{}$ & Shift &  $-0.5\,$ \\
		$\Sigma_{ij}^{}$ & Covariance & $\,\,\delta_{ij}^{}$ & $\rho + (1-\rho)\, \delta_{ij}^{}$ \\
		$K$ & Option strike & 1 \\
		$r$ & Interest rate & 0 \\
		\bottomrule
	\end{tabular}
	\vspace{2mm}
	\caption{Parameter settings in tests}
	\label{tbl:parameters}
	\vspace{-4mm}
\end{table}

\begin{figure}[!htbp]
	\centering
\begin{tikzpicture}
	\pgfplotsset{minor grid style={gray!10}}
	\pgfplotsset{major grid style={gray!30}}
	\begin{loglogaxis}[width=6.6cm, height=6.6cm, ylabel shift = 2mm, xlabel shift = 0mm, 
		xlabel={node count, $n$}, ylabel={error, $\left\| G - \tilde{G} \right\|_2^{}$}, 
		name=plot1, typeset ticklabels with strut,
		xmin=0.785, xmax=30, domain=1:1e3, ymin=1e-10, ymax=1e-2,
		legend pos=south west, xminorgrids, xmajorgrids, ymajorgrids,
		xtick style={draw=none}, ytick style={draw=none},
		xtick={1,10},
		xticklabels={1, 10},
		ytick={1e-10,1e-8,1e-6,1e-4,1e-2},
		yticklabels={$10^{-10}$,$10^{-8}$,$10^{-6}$,$10^{-4}$,$10^{-2}$,$10^{\,0~\,}$},
		legend cell align={left},
		]
		\addplot[gray!60, forget plot, smooth, domain=1:1e2] { 0.002 * exp(-1.42*x) };
		\addplot[gray!60, forget plot, smooth, domain=1:1e2] { 0.0035 * exp(-0.9*x) };
		\addplot[gray!60, forget plot, smooth, domain=1:1e2] { 0.007 * exp(-0.61*x) };
		%
		%
		\addplot[black, mark=+, mark options={fill=white}, only marks] table[x=n,y=error] {
			n	error
			1	0.003669406389932763
			3	0.0012335846159525687
			5	0.0003683009444743755
			7	0.00010552361129636964
			9	0.000030208186101344387
			11	8.026241758949207e-6
			15	6.730696730140577e-7
			21	2.2112211735979157e-8
		};
		\addlegendentry{~$\rho=0.5$}
		%
		\addplot[black, mark=asterisk, mark options={fill=white}, only marks] table[x=n,y=error] {
			n	error
			1	0.0014544054775031342
			3	0.00026408396224511836
			5	0.00004159973014979015
			7	6.4807983620905195e-6
			9	1.0289716236507137e-6
			11	1.6870945015683842e-7
			15	4.7415582500235856e-9
		};
		\addlegendentry{~$\rho=0.3$}
		%
		\addplot[black, mark=x, mark options={fill=white}, only marks] table[x=n,y=error] {
			n	error
			1	0.0004737853925527731
			3	0.000025618882594534998
			5	1.425310368952551e-6
			7	8.198882865943019e-8
			9	4.7226877478048794e-9
			11	3.605935145864696e-10
		};
		\addlegendentry{~$\rho=0.1$}
		\node[label={$\sim e^{-\alpha n}$}] at (17,3e-5){};
	\end{loglogaxis}
\end{tikzpicture}
\caption{Convergence of TT-X approximation error for correlated Gaussian distribution functions }
\label{fig:Gaussian}
\end{figure}
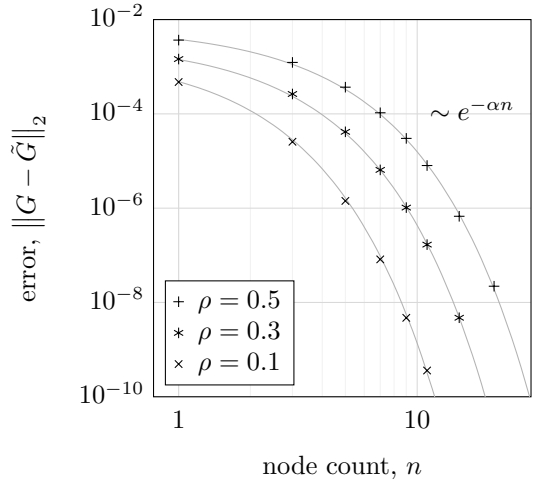

The accuracy of TT-X approximation for correlated Gaussian distribution functions ($G$ with $r=0$) is first presented in~\cref{fig:Gaussian}. Exponential convergence is observed for these analytic functions, where the error is measured as the root-mean-square residual using Monte Carlo. The choice of co-ordinate system is significant, but such exponential convergence is generally very efficient.
\vspace{2mm}

\begin{figure}[!htbp]
	\centering
\begin{tikzpicture}
	\pgfplotsset{minor grid style={gray!10}}
	\pgfplotsset{major grid style={gray!30}} 
	\begin{loglogaxis}[width=6.6cm, height=6.6cm, ylabel shift = 2mm, xlabel shift = 0mm, 
		xlabel={node count, $n$}, ylabel={error, $| \psi - \tilde{\psi} |$}, name=plot1, typeset ticklabels with strut,
		xmin=1, xmax=1e6, domain=1:1e12, ymin=1e-16, ymax=1,
		legend pos=south west, xminorgrids, xmajorgrids, ymajorgrids,
		xtick style={draw=none}, ytick style={draw=none},
		xtick={1,1e1,1e2,1e3,1e4,1e5,1e6,1e7,1e8,1e9},
		xticklabels={1, , $10^2$, , $10^4$, ,$10^6$, , ,$10^9$ },
		ytick={1e-15,1e-12,1e-9,1e-6,1e-3,1,1e3},
		yticklabels={$10^{-15}$,$10^{-12}$,$10^{-9}$,$10^{-6}$,$10^{-3}$,$10^{\,0~\,}$},
		legend cell align={left},
		]
		\addplot[gray!60, dashed, forget plot] { 0.65 * x^(-2) }; 
		\addplot[black] { 0.3081229594171603 * x^(-1/2) }; 
		\addlegendentry{Monte Carlo}
		\addplot[black, mark=*, mark options={fill=white}, only marks] table[x=n,y=e1] {
			n	t	e1	e2
			4	0.001556	0.017316347494803652	0.017316347494803652
			8	0.003709	0.008201988852291703	0.008309303180564998
			16	0.011105	0.0023293638350718737	0.00011271197213286177
			32	0.043157	0.00060453984174072	7.344503763773158e-7
			64	0.579443	0.00015649366784756258	9.455276128722545e-8
			128	1.472004	0.000040458665938158944	1.5671607234146734e-7
			256	3.952198	0.000010373799171736753	4.084484483968254e-8
			512	11.961318	2.637330962590445e-6	8.327190520462935e-9
			1024	44.858357	6.65861007442814e-7	1.6405301284105889e-9
			2048	180.695271	1.6719733109127333e-7	4.1337644418604214e-10
			4096	800.350693	4.175145829066196e-8	8.247560967511447e-11
			8192	3710.19411	1.0379504589153044e-8	4.3487757839244523e-10
			16384	17957.990426	2.82848336419228e-9	-1
		};
		\addlegendentry{TT-X}
		\addplot[black, mark=*, mark options={fill=black}, only marks] table[x=n,y=e2] {
			n	t	e1	e2
			16	0.011105	0.0023293638350718737	0.00011271197213286177
			32	0.043157	0.00060453984174072	7.344503763773158e-7
			64	0.579443	0.00015649366784756258	9.455276128722545e-8
			128	1.472004	0.000040458665938158944	1.5671607234146734e-7
			256	3.952198	0.000010373799171736753	4.084484483968254e-8
			512	11.961318	2.637330962590445e-6	8.327190520462935e-9
			1024	44.858357	6.65861007442814e-7	1.6405301284105889e-9
			2048	180.695271	1.6719733109127333e-7	4.1337644418604214e-10
			4096	800.350693	4.175145829066196e-8	8.247560967511447e-11
		};
		\addlegendentry{TT-X\,+\,Aitken}
		\node[label={$\sim n^{-\frac{1}{2}}$}] at (1e5,1e-3){};
		\node[label={$\sim n^{-2}$}] at (1e5,1e-11){}; 
	\end{loglogaxis}
	\begin{loglogaxis}[width=6.6cm, height=6.6cm, xlabel shift = 0mm, 
		xlabel={runtime, $\tau$/second}, ylabel={\phantom{,}}, yticklabels=\empty, 
		name=plot2, at=(plot1.right of south east), anchor=left of south west,  typeset ticklabels with strut,
		xmin=1e-4, xmax=1e6, domain=1e-4:1e8, ymin=1e-16, ymax=1,
		legend pos=south west, xminorgrids, xmajorgrids, ymajorgrids,
		xtick style={draw=none}, ytick style={draw=none},
		xtick={1e-4,1e-3,1e-2,1e-1,1e0,1e1,1e2,1e3,1e4,1e5,1e6,1e7,1e8},
		xticklabels={,$~~10^{-3}$, , , $1$, , ,$10^3$, , ,$10^6$ },
		ytick={1e-15,1e-12,1e-9,1e-6,1e-3,1,1e3},
		legend cell align={left},
		]
		\addplot[gray!60, dashed, forget plot, smooth, samples=100] { 0.00003 * x^(-1) };
		%
		%
		\addplot[black, smooth, samples=100] { 0.00031 * (x/2.8)^(-1/2) };
		\addlegendentry{Monte Carlo}
		\addplot[black, mark=*, mark options={fill=white}, only marks] table[x=t,y=e1] {
			n	t	e1	e2
			4	0.001556	0.017316347494803652	0.017316347494803652
			8	0.003709	0.008201988852291703	0.008309303180564998
			16	0.011105	0.0023293638350718737	0.00011271197213286177
			32	0.043157	0.00060453984174072	7.344503763773158e-7
			64	0.579443	0.00015649366784756258	9.455276128722545e-8
			128	1.472004	0.000040458665938158944	1.5671607234146734e-7
			256	3.952198	0.000010373799171736753	4.084484483968254e-8
			512	11.961318	2.637330962590445e-6	8.327190520462935e-9
			1024	44.858357	6.65861007442814e-7	1.6405301284105889e-9
			2048	180.695271	1.6719733109127333e-7	4.1337644418604214e-10
			4096	800.350693	4.175145829066196e-8	8.247560967511447e-11
			8192	3710.19411	1.0379504589153044e-8	4.3487757839244523e-10
			16384	17957.990426	2.82848336419228e-9	-1
		};
		\addlegendentry{TT-X}
		\addplot[black, mark=*, mark options={fill=black}, only marks] table[x=t,y=e2] {
			n	t	e1	e2
			16	0.011105	0.0023293638350718737	0.00011271197213286177
			32	0.043157	0.00060453984174072	7.344503763773158e-7
			64	0.579443	0.00015649366784756258	9.455276128722545e-8
			128	1.472004	0.000040458665938158944	1.5671607234146734e-7
			256	3.952198	0.000010373799171736753	4.084484483968254e-8
			512	11.961318	2.637330962590445e-6	8.327190520462935e-9
			1024	44.858357	6.65861007442814e-7	1.6405301284105889e-9
			2048	180.695271	1.6719733109127333e-7	4.1337644418604214e-10
			4096	800.350693	4.175145829066196e-8	8.247560967511447e-11
		};
		\addlegendentry{TT-X\,+\,Aitken}
		\node[label={$\sim \tau^{-\frac{1}{2}}$}] at (3e4,5e-6){}; 
		\node[label={$\sim \tau^{-1}$}] at (3e4,1.5e-11){};
	\end{loglogaxis}
\end{tikzpicture}
\vspace{-4mm}
\caption{Convergence of integration error }
\label{fig:ttxconvergence}
\end{figure}
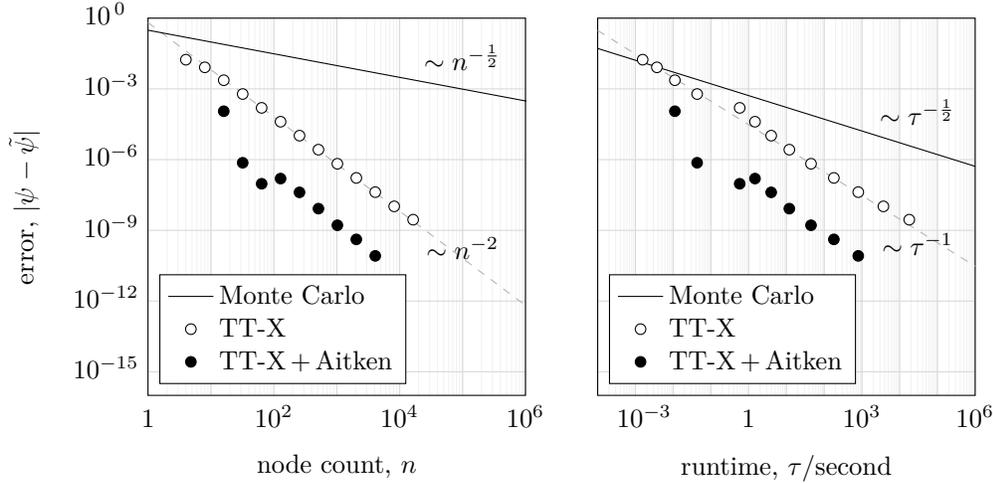

Not all TT-X representations converge exponentially though, and the integrand in~\cref{eq:psi} notably includes a non-analytic function. \cref{fig:ttxconvergence} demonstrates the associated integration value converging quadratically with respect to node count. Since the TT-X runtime roughly follows $\tau \propto n^2$ (due to matrix element evaluations) the respective convergence is linear in runtime.

The TT-X convergence rate still significantly outperforms Monte Carlo, and Aitken extrapolation (a general series acceleration technique) further suppresses errors by a few orders of magnitude. These convergence properties are problem-specific but they offer empirical proof of possible practical advantages.

The tests were run on a laptop without optimising the code efficiency, and so runtimes can be reduced. Furthermore, similar to Monte Carlo, many of the TT-X calculations can be computed in parallel. The (one-dimensional) TT-X integrations were performed analytically for this problem, but numerical integration techniques could otherwise be applied. On the error measure, the `exact' integration value~$\psi$ is deduced to machine precision using the Fourier-TT network representation (described next).

\vspace{3mm}\noindent
\textbf{Fourier Tensor Train (Fourier-TT) network}

\noindent
Other constructive methods for forming tensor networks can be based on discretisation of integral transforms. However, such approaches are problem-specific. For the function $\phi$ in~\cref{eq:phi}, it is convenient to consider the Fourier series of an effective one-dimensional problem:
\begin{eqnarray}
	\max [0, S-K] &~=\,& \left( \frac{d-1}{ d } \right) S - \! \sum_{m=1}^{\infty} \frac{ 2 Kd }{m^2 \, \pi^2} 
	\, \sin \Big[ \frac{ m \pi }{ d } \Big] \, \sin \Big[ \frac{ m \pi S }{ Kd } \Big]
	\hspace{3mm}
	S \in [0, Kd]
	\hspace{5mm}\phantom{.}
\end{eqnarray}
On substituting $S = \sum_i^d \, \omega_i^{} \, e^{x_i^{}}$, the multivariate function with respect to $\{ x_i^{} \}$ is recovered and the series representation is valid in the domain $\{ \omega_i^{}\, e^{x_i^{}} \} \in [0,K]^d$. Given that the relevant Green's function can be efficiently expressed in a separable form using its TT-X representation, this domain encloses the `difficult' part of the integration problem. Outside the domain, $\phi$ is already separable (for non-negative $\omega_i^{}$).

\begin{figure}[!htbp]
\centering
\begin{tikzpicture}
\pgfplotsset{minor grid style={gray!10}}
\pgfplotsset{major grid style={gray!30}} 
\begin{loglogaxis}[width=6.6cm, height=6.6cm, ylabel shift = 2mm, xlabel shift = 0mm, 
	xlabel={Fourier-TT series terms, $n$}, ylabel={error, $| \psi - \tilde{\psi} |$}, name=plot1, typeset ticklabels with strut,
	xmin=1, xmax=1e4, domain=1:1e12, ymin=1e-15, ymax=1,
	legend pos=south west, xminorgrids, xmajorgrids, ymajorgrids,
	xtick style={draw=none}, ytick style={draw=none},
	xtick={1,1e1,1e2,1e3,1e4,1e5,1e6,1e7,1e8,1e9},
	xticklabels={1, , $10^2$, , $10^4$, ,$10^6$, , ,$10^9$ },
	ytick={1e-15,1e-12,1e-9,1e-6,1e-3,1,1e3},
	yticklabels={$10^{-15}$,$10^{-12}$,$10^{-9}$,$10^{-6}$,$10^{-3}$,$10^{\,0~\,}$},
	legend cell align={left},
	]
	\addplot[gray!60, forget plot, smooth, domain=100:1e4] { 0.15 * exp(-0.011*x) };
	\addplot[gray!60, forget plot, smooth, domain=10:1e3] { 0.05 * exp(-0.14*x) };
	%
	%
	\addplot[black, mark=x, mark options={fill=white}, only marks] table[x=n,y=error] {
		n	error
		100	0.054112886885206083867689
		200	0.014581851566824181678194
		300	0.004513184856529696312253
		400	0.00132100685091316919758
		500	0.000305477723202235702272
		600	0.000023881501929959222694
		700	0.000026723004261446307336
		800	0.000020126228937498258738
		900	8.571223705598221481e-6
		1000	2.32383971616774946e-6
		1100	2.08634037728120128e-7
		1200	1.58749819426641648e-7
		1300	9.2869326590043436e-8
		1400	2.12851627918520182e-8
		1500	1.5842775991644277e-9
		1600	2.7278294487712332e-9
		1700	7.100673436978262e-10
		1800	6.67504007180743e-11
		1900	8.843044702587619e-11
		2000	1.447189580972726e-11
		2100	5.67452804701451e-12
		2200	2.512159439370791e-12
		2300	1.44713900502373e-13
		2400	2.71023905413262e-13
		2500	1.87922269902181e-14
		2600	2.53075146553322e-14
		2700	3.9601947875963e-15
		2800	2.32768366740731e-15
		2900	5.0483830298076e-16
		3000	2.27867231791607e-16
		3100	5.3901479606242e-17
		3200	2.45038474716156e-17
		3300	5.0076394424268e-18
		3400	2.8409100458818e-18
		3500	3.6618605299777e-19
		3600	3.3722404424548e-19
		3700	9.170116409373e-21
		3800	3.8405227055517e-20
		3900	3.6447727059336e-21
		4000	3.8538262385595e-21
	};
	%
	\addplot[black, mark=asterisk, mark options={fill=white}, only marks] table[x=n,y=error] {
		n	error
		10	0.010440452415119326723882475
		20	0.00158043801291385327715217758
		30	0.000346823206456692808289478304
		40	0.00006458270131400928894228591
		50	7.925018702204823086721043e-6
		60	5.63029237879059438872938259e-6
		70	1.104347357255418516665862218e-6
		80	7.3406632152066786068808268e-8
		90	1.480296761359231124280846288e-7
		100	6.533836029670690233225093014e-8
		110	1.68027524893315674770624205e-8
		120	1.076017349752252244339375239e-9
		130	2.02461519396438018346306248e-9
		140	1.541807333975180656765188271e-9
		150	7.67679360245433899588008091e-10
		160	2.94483758871948893856583308e-10
		170	7.5288744065745178572154867e-11
		180	2.057937882551786431944595e-12
		190	1.6789591922599525903497332e-11
		200	1.46242708771442843365568862e-11
		210	9.52159219299585884510310863e-12
		220	5.08788442428144640906214427e-12
		230	2.35534217747700880430857467e-12
		240	9.3826814661200675314456003e-13
		250	2.4345146779070347756104473e-13
		260	4.841458620103396076858542e-15
		270	9.300734910984313837115143e-14
		280	8.8240324475396570933191391e-14
		290	7.033863330112872159767973e-14
		300	4.55089408664342163025687003e-14
		310	2.81249502994121074461361128e-14
		320	1.56110890421023601666074494e-14
		330	8.0264134456426594103404628e-15
		340	3.759557479491505440662359e-15
		350	1.3436458505574047773606354e-15
		360	3.1146987731204637083647979e-16
		370	2.3444719668849493074817834e-16
		380	3.2636859207674715906977053e-16
		390	3.57949716570376236488005003e-16
		400	2.73434562526162133606662256e-16
		410	2.17356463265470074651749028e-16
		420	1.4494678956926048679080749e-16
		430	1.00789613957127651563650271e-16
		440	6.20800026170400117440771742e-17
		450	3.88274801610108985496726166e-17
		460	2.2128257478066929113902934e-17
		470	1.2011096975131339236783463e-17
		480	6.0046617712778675324800292e-18
		490	2.23570616373458704079429592e-18
		500	5.7879157795885632693212137e-19
	};
	\node[label={$d=100$}] at (8e2,5e-3){};
	\node[label={$d=10$}] at (1e1,5e-3){};
	\node[label={$\sim e^{-\alpha n}$}] at (4e1,1e-15){};
\end{loglogaxis}
\end{tikzpicture}
\caption{Convergence of Fourier-TT integration error }
\label{fig:fouriertt}
\end{figure}
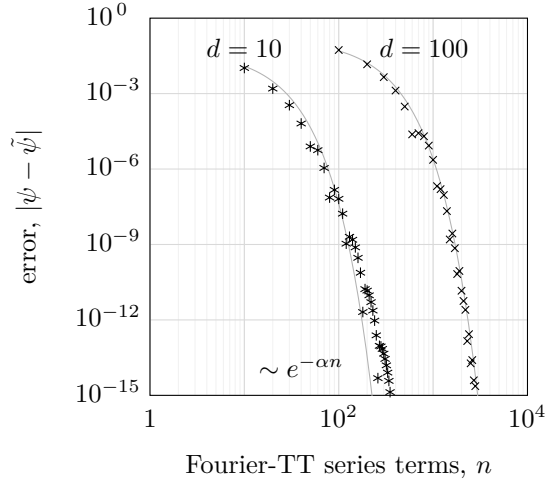

In order to then express the Fourier series in a separable form, a tensor train representation of the sine function can be applied. An exact representation is related to the product of rotation matrices, where edge matrices can be collapsed to the relevant row/column in order to reproduce the target sine function.
\begin{eqnarray}
	\mathbf{R} \Bigg[ \sum_{i}^{d} \theta_i^{} \Bigg] &~=~&
	\mathbf{R} [\theta_{1}^{}] \cdots \mathbf{R} [\theta_{d}^{}]
	\hspace{4mm} \text{where} \hspace{2mm}
	\mathbf{R} [\theta] ~=\,
	\begin{pmatrix}
		\cos \theta & -\sin \theta \\
		\sin \theta& \phantom{-} \cos \theta
	\end{pmatrix}
\end{eqnarray}
This construction is designated as a Fourier tensor train (Fourier-TT) representation in this article, and \cref{fig:fouriertt} demonstrates that exponential convergence is roughly achieved for the integration problem. In this test, numerical integration is applied for the relevant (one-dimensional) integrals.

\vfill\pagebreak

\cref{tbl:convergence} summarises the runtime orders of magnitude for different methods. For precise calculations, the tensor networks offer significant advantages with respect to Monte Carlo, given superior convergence properties. However, if low precision results are acceptable, Monte Carlo may still be preferred.
\vspace{3mm}
\begin{table}[!htbp]
	\centering
	\begin{tabular}{lllll} 
		\toprule
		&& \multicolumn{3}{c}{Convergence time} \\
		&  & $\epsilon < 10^{-3}$ & $\epsilon < 10^{-6}$ & $\epsilon < 10^{-9}$ \\ \midrule
		\multirow{4}{*}{$d=10$} 
		&Monte Carlo & 0.1\,sec	&1~day	& $10^{4}$~years \\
		&TT-X & 0.01\,sec	& 10~sec	& 10~hours \\
		&TT-X\,+\,Aitken & 0.01\,sec	& 0.1\,sec	& 100~sec \\
		&Fourier-TT & 0.1\,sec	& 0.1\,sec	& 0.1~sec \\ \midrule
		\multirow{4}{*}{$d=100$} 
		&Monte Carlo & 0.1\,sec	&  1~day	&  $10^{4}$~years  \\
		&TT-X & 1~sec	&  1~hour	&  10~days  \\
		&TT-X\,+\,Aitken & 1~sec	&  10~sec	&  10~min  \\
		&Fourier-TT & 1~sec	&  5~sec	&  10~sec  \\
		\bottomrule
	\end{tabular}
	\vspace{2mm}
	\caption{Convergence times for integration error $\epsilon = | \psi - \tilde{\psi} |$ } 
\label{tbl:convergence}
\end{table}

\vspace{-2mm}
\vspace{3mm}\noindent
\textbf{Summary}

\noindent
This article detailed regression-free methods for forming tensor network representations of functions, and demonstrated practical benefits when used for high-dimensional integration. Although the convergence properties are problem-specific, tensor networks were shown to converge faster than Monte Carlo for a typical problem in finance: basket option valuation. The tensor train cross representation can be used for any integrand, and so this method is applicable for any problem that can be expressed as an integral. Generally, tensor networks offer opportunities for efficient integration.

\vfill\pagebreak
\noindent\textbf{References}\vspace{1mm}
\renewcommand{\refname}{}
\renewcommand{\bibsection}{}
\setlength{\bibsep}{1mm}

\end{document}